\documentclass[12pt]{article}
\usepackage{amssymb}
\usepackage{latexsym,bm}
\usepackage{graphicx}
\usepackage{amsmath}
\usepackage{mathrsfs}
\usepackage{epstopdf}

\date{}
\newcounter{mathitem}
  {\begin{list}{{$(\roman{mathitem})$}}{
   \setcounter{mathitem}{0}
   \usecounter{mathitem}
   \setlength{\topsep}{0pt plus 2pt minus 0pt}
   \setlength{\parskip}{0pt plus 2pt minus 0pt}
   \setlength{\partopsep}{0pt plus 2pt minus 0pt}
   \setlength{\parsep}{0pt plus 2pt minus 0pt}
   \setlength{\leftmargin}{35pt}
   \setlength{\itemsep}{0pt plus 2pt minus 0pt}}}
  {\end{list}}

\begin{document}
\title{The minimum size of a $3$-connected locally nonforesty graph
\footnote{E-mail addresses: {\tt lichengli0130@126.com} (C. Li), {\tt tyr2290@163.com} (Y. Tang), {\tt zhan@math.ecnu.edu.cn} (X. Zhan).}}
\author{\hskip -10mm Chengli Li, Yurui Tang and Xingzhi Zhan\thanks{Corresponding author}\\
{\hskip -10mm \small Department of Mathematics, East China Normal University, Shanghai 200241, China}}\maketitle
\begin{abstract}
A local subgraph of a graph is the subgraph induced by the neighborhood of a vertex. Thus a graph of order $n$ has $n$ local subgraphs. A graph $G$ is called locally nonforesty
 if every local subgraph of $G$ contains a cycle. Recently, in studying forest cuts of a graph, Chernyshev, Rauch and Rautenbach posed the conjecture that if $n$ and $m$ are the order and size of a $3$-connected locally nonforesty graph respectively, then $m\ge 7(n-1)/3.$ We solve this problem by determining the minimum size of a $3$-connected locally nonforesty graph of order $n.$ It turns out that the conjecture does not hold.
\end{abstract}


\vskip 8mm

\section{Introduction}

We consider finite simple graphs and use standard terminology and notation from [2] and [9]. The {\it order} of a graph is its number of vertices, and the
{\it size} its number of edges.

{\bf Definition 1.} A {\it local subgraph} of a graph $G$ is the subgraph of $G$ induced by the neighborhood of a vertex.

Thus a graph of order $n$ has $n$ local subgraphs. It has been a traditional topic to deduce properties of a graph by its local subgraphs. A graph $G$ is called {\it locally connected} if all local subgraphs of $G$ are connected. {\it Locally hamiltonian} graphs are defined similarly [5]. A well-known theorem of Oberly and Sumner [7] states that every connected, locally connected claw-free graph is hamiltonian. Further connectivity conditions on such graphs force them to be hamiltonian-connected ([1] and [3]). An interesting unsolved conjecture of Ryj\'{a}\v{c}ek (see [8]) asserts that every locally connected graph is weakly pancyclic.

{\bf Definition 2.} A graph $G$ is called {\it locally foresty} if every local subgraph of $G$ is a forest. $G$ is called {\it locally nonforesty} if every local subgraph
of $G$ contains a cycle.

Recently, in studying forest cuts of a graph, Chernyshev, Rauch and Rautenbach [4, p.8] posed the following

{\bf Conjecture 1.} {\it If $n$ and $m$ are the order and size of a $3$-connected locally nonforesty graph respectively, then $m\ge 7(n-1)/3.$}

In this paper we determine the minimum size of a $3$-connected locally nonforesty graph of order $n.$ It turns out that Conjecture 1 does not hold.

The dual problem is to ask: What is the maximum size of a locally foresty  graph of order $n?$ Note that a graph is locally foresty if and only if it contains
no wheels. This problem has been solved by Moon [6].

We denote by $V(G)$ the vertex set of a graph $G,$ and denote by $e(G)$ the size of $G.$ The neighborhood of a vertex $x$ is denoted by $N(x)$ or $N_G(x),$ and the
 closed neighborhood of $x$ is $N[x]\triangleq N(x)\cup \{x\}.$ The degree of $x$ is denoted by ${\rm deg}(x).$ We denote by $\delta (G)$ and $\Delta(G)$ the minimum degree and maximum degree of $G,$ respectively.  For a vertex subset $S\subseteq V(G),$ we use $G[S]$ to denote the subgraph of $G$ induced by $S,$ and use $N(S)$ to denote the neighborhood
 of $S;$ i.e., $N(S)=\{y\in V(G)\setminus S \,|\, y\,\,{\rm has}\,\,{\rm a}\,\,{\rm neighbor}\,\,{\rm in}\,\,S\}.$ For $x\in V(G)$ and $S\subseteq V(G),$ $N_S(x)\triangleq N(x)\cap S$ and the degree of $x$ in $S$ is ${\rm deg}_S(x)\triangleq |N_S(x)|.$ Given two disjoint vertex subsets $S$ and $T$ of $G,$ we denote by $[S, T]$ the set of edges having one endpoint in $S$ and the other in $T.$  The degree of $S$ is ${\rm deg}(S)\triangleq |[S, \overline{S}]|,$ where $\overline{S}=V(G)\setminus S.$ We denote by $C_n,$ $P_n$ and $K_n$ the cycle of order $n,$ the path of order $n$ and the complete graph of order $n,$ respectively. $\overline{G}$ denotes the complement of a graph $G.$ For two graphs $G$ and $H,$ $G\vee H$ denotes the {\it join} of $G$ and $H,$ which is obtained from the disjoint union $G+H$ by adding edges joining every vertex of $G$ to every vertex of $H.$

For graphs we will use equality up to isomorphism, so $G=H$ means that $G$ and $H$ are isomorphic.

\section{The minimum size of a $3$-connected locally nonforesty graph}

In this section we determine the minimum size of a $3$-connected locally nonforesty graph of order $n.$

{\bf Theorem 1.} {\it Given an integer $n\ge 8,$ define
$$
f(n)=\begin{cases}
2n-\lfloor n/8\rfloor\quad\quad\quad\,{\rm if}\,\,\, n\equiv 0,\,4,\,7\,(\,{\rm mod}\,\, 8),\\
2n+1-\lfloor n/8\rfloor \quad\,\,\,{\rm otherwise.}
\end{cases}
$$

Then the minimum size of a $3$-connected locally nonforesty graph of order $n$ is $f(n).$
}

{\bf Proof.} Let $G$ be a $3$-connected locally nonforesty graph of order $n.$ We first prove that $e(G)\ge f(n).$
Let $n=8k+r$ with $0\le r\le 7.$ Then
$$
f(n)=\begin{cases}
15k+2r\quad\quad\quad\,{\rm if}\,\,\, r=0,\,4,\,7,\\
      15k+2r+1 \quad\,\,\,{\rm otherwise.}
\end{cases}
$$
For a vertex $v$ of $G,$ we denote by $L(v)$ the subgraph of $G$ induced by $N(v).$

We will repeatedly use the degree-sum formula for the size of a graph ([2, p.7] or [9, p.35]. Since $G$ is $3$-connected, $\delta(G)\ge 3.$ If $\delta(G)\ge 4,$ we have
$e(G)\ge 2n\ge f(n).$ Next we suppose $\delta(G)=3.$ Denote
$$
S=\{v\in V(G)|\,{\rm deg}(v)=3\},\quad T=N(S),\quad W=V(G)\setminus (S\cup T).
$$
We assert that $S$ is an independent set. To the contrary, suppose that $S$ contains two adjacent vertices $u_1$ and $u_2.$ Since $L(u_1)$ contains a cycle,
$G[N[u_1]]=K_4$ which contains $u_2.$ Then $N(u_1)\setminus \{u_2\}$ would be a vertex cut of $G$ of cardinality $2,$ contradicting our assumption that $G$ is $3$-connected.

Let $T_i=\{v\in T|\, {\rm deg}_S(v)=i\}$ where $1\le i\le s\triangleq |S|.$ For any vertex $u\in T_i,$ we have ${\rm deg}(u)\ge i+2,$ which can be deduced from the facts
that $N_S(u)$ is an independent set and $L(u)$ contains a cycle. We distinguish two cases.

{\bf Case 1.} There exists an $i$ such that $T_i$ contains a vertex $u$ with ${\rm deg}(u)=i+2.$

Let $A=N_S(u)$ where $|A|=i$ and $N(u)\setminus S=\{v, w\}.$ For each $p\in A,$ $G[N[p]]=K_4,$ implying that $N(p)\subseteq N[u]=A\cup\{u,v,w\}.$
Since $A$ is an independent set, we deduce that $N(p)=\{u,v,w\}$ for any $p\in A.$ If there is a vertex not in $A\cup \{u,v,w\},$ then $\{v,w\}$ would be a vertex cut of $G,$ contradicting the condition that $G$ is $3$-connected. Hence $V(G)=A\cup \{u,v,w\},$ implying that $G=K_3\vee \overline{K_i}.$ It follows that $e(G)=3n-6>f(n),$ where 
we have used the assumption that $n\ge 8.$

{\bf Case 2.} For every $i$ and any vertex $u\in T_i,$ ${\rm deg}(u)\ge i+3.$

Recall that we have denoted $s=|S|.$ Note that $\sum_{i=1}^s |T_i|=|T|$ and
$$
\sum_{i=1}^s i\cdot |T_i|=\sum_{v\in T}{\rm deg}_S(v)=|[S,\, T]|=\sum_{u\in S}{\rm deg}_T(u)=3s.
$$
Observe that every vertex in $T\cup W$ has degree at least $4.$ We have
$$
e(G)\ge (3s+4(n-s))/2=2n-s/2
$$
and
\begin{align*}
e(G)&\ge (3s+\sum_{i=1}^{s} (i+3)|T_i|+4|W|)/2\\
    &=(6s+3|T|+4|W|)/2\\
    &\ge 3(n+s)/2.
\end{align*}
Denote $\varphi(s)={\rm max}\{\lceil 2n-s/2\rceil,\,\lceil 3(n+s)/2\rceil \}.$ Then
$$
e(G)\ge \mathop{\rm min}\limits_{1\le s\le n}\varphi(s). \eqno (1)
$$
According to the value of the remainder $r$ in $n=8k+r$ with $0\le r\le 7,$ we distinguish eight subcases.

{\bf Subcase 2.1.} $n=8k.$

By inequality (1), we have $e(G)\ge \mathop{\rm min}\limits_{1\le s\le n}\varphi(s)=\varphi(2k)=15k=f(n).$

{\bf Subcase 2.2.} $n=8k+1.$

In this case, $\mathop{\rm min}\limits_{1\le s\le n}\varphi(s)=15k+2$ and the minimum value $15k+2$ is attained uniquely at $s=2k.$ By (1) we have
$e(G)\ge 15k+2$ and equality holds only if $s=2k.$ Now we exclude the possibility that $e(G)=15k+2,$ and then it follows that $e(G)\ge 15k+3=f(n).$

To the contrary, suppose $e(G)=15k+2.$ Then $s=2k.$ It follows that the degree sequence of $G$ is
$(\operatorname*{\underbrace{3,\cdots, 3}}\limits_{2k},\operatorname*{\underbrace{4,\cdots, 4}}\limits_{6k+1}).$ Let $S=\{a_1,\ldots,a_{2k}\},$
the set of vertices of degree $3.$ Denote $A_i=N[a_i],$ the closed neighborhood of $a_i.$ Then $G[A_i]=K_4$ for each $i=1,\ldots,2k.$

{\bf Claim 1.} $A_i\cap A_j=\emptyset$ for $i\neq j.$

To the contrary, suppose  $A_i\cap A_j\neq\emptyset.$ Recall that $S=\{a_1,\ldots,a_{2k}\}$ is an independent set.  We have $|A_i\cap A_j|\le 3.$
If $|A_i\cap A_j|=1$ or $2,$ then  $A_i\cap A_j$ contains a vertex of degree at least $5,$ a contradiction. If $|A_i\cap A_j|=3$ then $|V(G)|=5,$ contradicting our
assumption that $|V(G)|\ge 8.$

Let $V(G)\setminus (\bigcup_{i=1}^{2k}A_i)=\{x\}.$

{\bf Claim 2.} $|N(x)\cap A_i|\le 1$ for every $i$ with $1\le i\le 2k.$

Note that $|A_i|=4$ and $a_i$ and $x$ are nonadjacent for each $i.$ Thus $|N(x)\cap A_i|\le 3.$ Using the conditions $\Delta(G)=4$ and ${\rm deg}(a_i)=3,$ we deduce that
if $|N(x)\cap A_i|=3$ then $x$ is a cut-vertex of $G,$ and if $|N(x)\cap A_i|=2,$ with $N(a_i)\setminus N(x)=\{u\},$ $\{x,\,u\}$ is a vertex cut of $G,$ contradicting
the condition that $G$ is $3$-connected. This shows that $|N(x)\cap A_i|\le 1.$

Claim 2 implies that $N(x)$ is an independent set and hence $L(x)$ does not contain a cycle, a contradiction. This shows that $e(G)\neq 15k+2.$

In the following six subcases, Subcases 2.5 and 2.8 are easily treated as Subcase 2.1, while the remaining four subcases are treated by using ideas similar to those
in Subcase 2.2, but there are new subtle differences.

{\bf Subcase 2.3.} $n=8k+2.$

$f(n)=15k+5.$ Now $\mathop{\rm min}\limits_{1\le s\le n}\varphi(s)=15k+4$ and the minimum value $15k+4$ is attained uniquely at $s=2k.$ The proof in this case is very similar to
the one in Subcase 2.2, and we omit the details.

{\bf Subcase 2.4.} $n=8k+3.$

$f(n)=15k+7.$ We have $\mathop{\rm min}\limits_{1\le s\le n}\varphi(s)=15k+6$ and the minimum value $15k+6$ is attained at $s=2k$ or $s=2k+1.$
By (1) we have $e(G)\ge 15k+6$ and equality holds only if $s=2k$ or $s=2k+1.$ It suffices to exclude the possibility that $e(G)=15k+6.$ To the contrary, suppose $e(G)=15k+6.$
Then $s=2k$ or $s=2k+1.$ The case $s=2k$ is similar to Subcase 2.2 and we omit the details. Now suppose $s=2k+1.$ Then the degree sequence of $G$ is
$(\operatorname*{\underbrace{3,\cdots, 3}}\limits_{2k+1},\operatorname*{\underbrace{4,\cdots, 4}}\limits_{6k+1},5).$ We define $a_i$ and $A_i,$ $1\le i\le 2k+1$ as above.
Then $G[A_i]=K_4$ for each $i.$

{\bf Claim 3.} $A_i\cap A_j=\emptyset$ for $i\neq j.$

To the contrary, suppose  $A_i\cap A_j\neq\emptyset.$ Note that $\Delta(G)=5.$ If $|A_i\cap A_j|=1,$ then the vertex in $A_i\cap A_j$ has degree at least $6,$ a contradiction.
If $|A_i\cap A_j|=2,$ then the two vertices in $A_i\cap A_j$ have degree at least $5,$ contradicting the degree sequence of $G.$ If $|A_i\cap A_j|=3,$ then one vertex in $A_i$
has degree $5$ and that vertex is a cut-vertex of $G,$ a contradiction.

Claim 3 yields that $|V(G)|\ge 4(2k+1)=8k+4,$ contradicting $n=8k+3.$ This shows that $e(G)\neq 15k+6.$

{\bf Subcase 2.5.} $n=8k+4.$

We have $e(G)\ge \mathop{\rm min}\limits_{1\le s\le n}\varphi(s)=\varphi(2k)=15k+8=f(n).$

{\bf Subcase 2.6.} $n=8k+5.$

$f(n)=15k+11.$ We have $\mathop{\rm min}\limits_{1\le s\le n}\varphi(s)=15k+10$ and the minimum value $15k+10$ is attained at $s=2k$ or $s=2k+1.$ It suffices to exclude the possibility that $e(G)=15k+10.$ To the contrary, suppose $e(G)=15k+10.$ Then $s=2k$ or $s=2k+1.$

We first consider the case $s=2k.$ The degree sequence of $G$ is
$(\operatorname*{\underbrace{3,\cdots, 3}}\limits_{2k},\operatorname*{\underbrace{4,\cdots, 4}}\limits_{6k+5}).$

We define $a_i$ and $A_i,$ $1\le i\le 2k$ as above. As in Subcase 2.2 we have the following

{\bf Claim 4.} $A_i\cap A_j=\emptyset$ for $i\neq j.$

Denote $R=V(G)\setminus (\bigcup_{i=1}^{2k}A_i).$ Since $n=8k+5,$ Claim 4 implies $|R|=5.$ As in Subcase 2.2 we can prove the following

{\bf Claim 5.} $|N(y)\cap A_i|\le 1$ for any $y\in R$ and any $i$ with $1\le i\le 2k.$

By Claim 5, the subgraph $H\triangleq G[R]$ is a locally nonforesty graph. It is easy to determine the minimum size of such a graph of order $5,$
and we have $e(H)\ge 9.$ We estimate the degree of $R$ as follows:
\begin{align*}
{\rm deg}(R)=\sum_{v\in R}({\rm deg}(v)-{\rm deg}_R(v))=&\sum_{v\in R}{\rm deg}(v)-\sum_{v\in R}{\rm deg}_R(v)\\
            &=20-2e(H)\le 2,
\end{align*}
which implies that the connectivity of $G$ is at most $2,$ contradicting our assumption that $G$ is $3$-connected.

Now we analyze the case $s=2k+1.$ The degree sequence of $G$ is
$(\operatorname*{\underbrace{3,\cdots, 3}}\limits_{2k+1},\operatorname*{\underbrace{4,\cdots, 4}}\limits_{6k+3},5).$

We define $a_i$ and $A_i,$ $1\le i\le 2k+1$ as above. As in Subcase 2.4 we have the following

{\bf Claim 6.} $A_i\cap A_j=\emptyset$ for $i\neq j.$

Let $V(G)\setminus (\bigcup_{i=1}^{2k+1}A_i)=\{p\}.$ Let $u$ be the vertex of degree $5.$

{\bf Claim 7.} $|N(p)\cap A_i|\le 1$ for every $i$ with $1\le i\le 2k+1.$

We fix an $i.$ If $u\notin A_i,$ the proof of Claim 7 is the same as in Subcase 2.2. Next suppose $u\in A_i.$ If $|N(p)\cap A_i|=3,$ then $\{p, u\}$ is a vertex cut of $G,$
a contradiction. Suppose $|N(p)\cap A_i|=2.$  If $p$ and $u$ are nonadjacent, then $\{p, u\}$ is a vertex cut of $G,$ a contradiction. If $p$ and $u$ are adjacent, using the conditions that ${\rm deg}(u)=5$ and ${\rm deg}(p)=4$ we deduce that $L(p)$ cannot contain a cycle, a contradiction.

Claim 7 implies that $N(p)$ is an independent set. Hence $L(p)$ is an empty graph, a contradiction.

{\bf Subcase 2.7.} $n=8k+6.$

$f(n)=15k+13.$ We have $\mathop{\rm min}\limits_{1\le s\le n}\varphi(s)=15k+12$ and the minimum value $15k+12$ is attained at $s=2k,$ $s=2k+1$ or $s=2k+2.$ It suffices to
exclude the possibility that $e(G)=15k+12.$ To the contrary, suppose $e(G)=15k+12.$ Then $s=2k,$ $s=2k+1$ or $s=2k+2.$

(a) Assume $s=2k.$

Now the degree sequence of $G$ is
$(\operatorname*{\underbrace{3,\cdots, 3}}\limits_{2k},\operatorname*{\underbrace{4,\cdots, 4}}\limits_{6k+6}).$
It is easy to verify that the minimum size of a locally nonforesty graph of order $6$ is $11.$ Arguing as in Subcase 2.6 we deduce that the connectivity of $G$ is at most $2,$
a contradiction.

(b) Assume $s=2k+1.$

The degree sequence of $G$ is
$(\operatorname*{\underbrace{3,\cdots, 3}}\limits_{2k+1},\operatorname*{\underbrace{4,\cdots, 4}}\limits_{6k+4},5).$ Arguing as in Subcase 2.6 we can find a vertex whose
neighborhood induces a forest, a contradiction.

(c) Assume $s=2k+2.$

The degree sequence of $G$ is either
$(\operatorname*{\underbrace{3,\cdots, 3}}\limits_{2k+2},\operatorname*{\underbrace{4,\cdots, 4}}\limits_{6k+2},5,5)$ or
$(\operatorname*{\underbrace{3,\cdots, 3}}\limits_{2k+2},\operatorname*{\underbrace{4,\cdots, 4}}\limits_{6k+3},6).$ We define $a_i$ and $A_i,$ $1\le i\le 2k+2$ as before.
Since $n=8k+6$ and $\sum_{i=1}^{2k+2}|A_i|=8k+8,$ there exist $j\neq q$ such that $A_j\cap A_q\neq\emptyset.$ Since the set $S=\{a_1,\ldots,a_{2k+2}\}$ is an independent set,
we have $|A_j\cap A_q|\le 3.$

Suppose $|A_j\cap A_q|=3;$ i.e., $N(a_j)=N(a_q).$ Let $F$ be the set of vertices in $N(a_j)$ that have degree greater than $4.$ In view of the degree sequence of $G$ and
using the fact that $n\ge 8,$ we deduce that $1\le |F|\le 2.$ But then $F$ is a vertex cut of $G$ of cardinality at most $2,$ contradicting the condition that $G$ is $3$-connected.

Suppose $|A_j\cap A_q|=2.$ Recall that the symmetric difference of two sets $A$ and $B,$ denoted $A\triangledown B,$ is the set $(A\cup B)\setminus (A\cap B).$
Since $n\ge 8,$ the set $N(a_j)\triangledown N(a_q)$ is a vertex cut of $G$ of cardinality $2,$ a contradiction.

Suppose $|A_j\cap A_q|=1.$ Since $n=8k+6$ and $\sum_{i=1}^{2k+2}|A_i|=8k+8,$ besides $(j, q),$ there exists another index pair $(h, r)$ (possibly $|\{j, q\}\cap \{h, r\}|=1$)
such that $|A_h\cap A_r|=1.$ Let $A_j\cap A_q=\{v_1\}$ and $A_h\cap A_r=\{v_2\}.$ Then ${\rm deg}(v_1)\ge 6$ and ${\rm deg}(v_2)\ge 6.$
If $v_1=v_2,$ then $v_1$ has degree at least $9,$ contradicting the degree sequence of $G;$ if $v_1\neq v_2,$ then $G$ has at least two vertices of degree at least $6,$
contradicting again the degree sequence of $G.$

{\bf Subcase 2.8.} $n=8k+7.$

We have $e(G)\ge \mathop{\rm min}\limits_{1\le s\le n}\varphi(s)=\varphi(2k)=15k+14=f(n).$

Next for every integer $n\ge 8,$ we construct a $3$-connected graph $G_n$ of order $n$ and size $f(n)$ each of whose local subgraph contains a cycle. Thus $f(n)$ is indeed the
minimum size.

The graphs $G_8, G_9, G_{10}$ and $G_{11}$ are depicted in Figure 1, and the  graphs $G_{12}, G_{13}, G_{14}$ and $G_{15}$ are depicted in Figure 2.
\begin{figure}[h]
\centering
\includegraphics[width=0.6\textwidth]{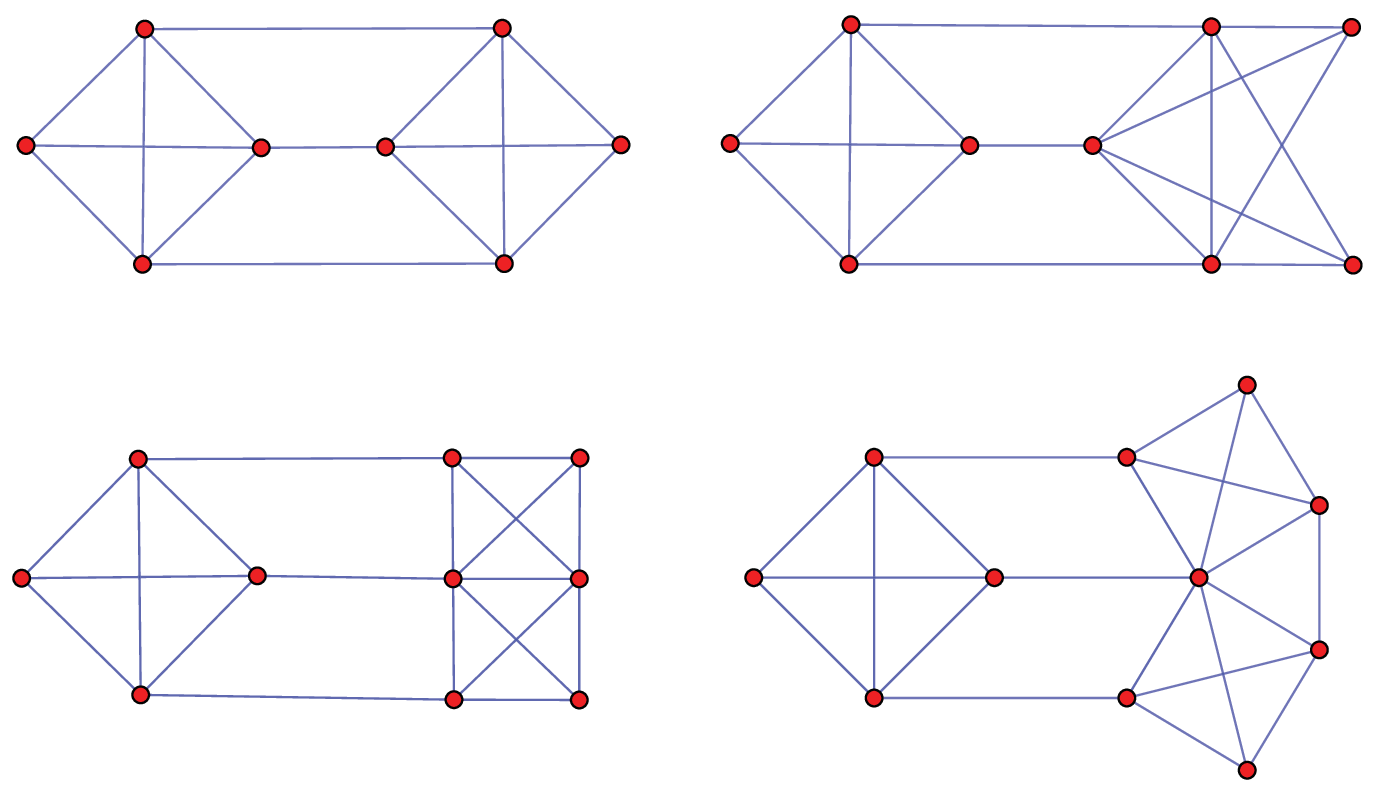}
\caption{$G_8, G_9, G_{10}$ and $G_{11}$}
\end{figure}

\begin{figure}[h]
\centering
\includegraphics[width=0.5\textwidth]{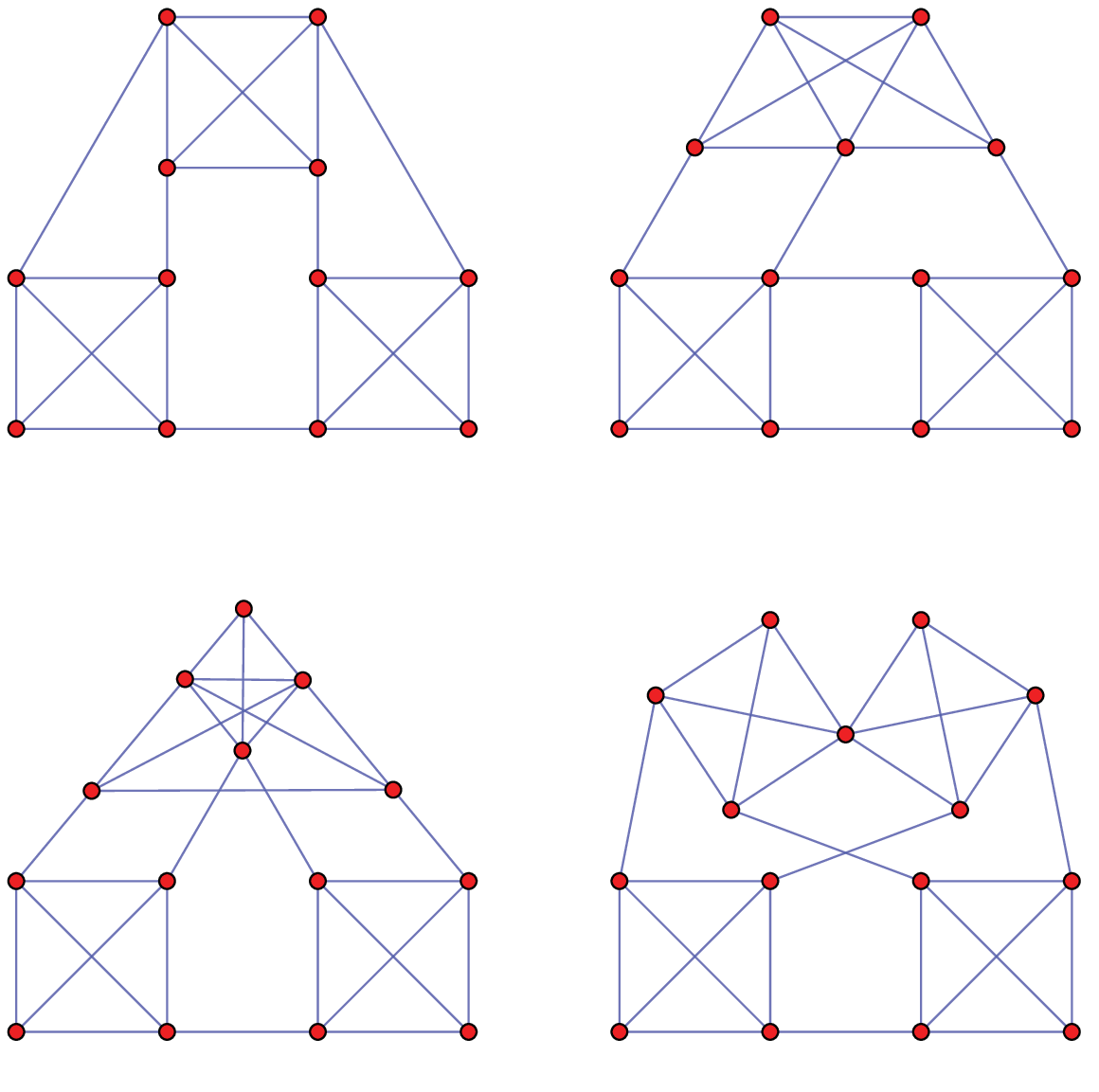}
\caption{$G_{12}, G_{13}, G_{14}$ and $G_{15}$}
\end{figure}

Now suppose $n\ge 16.$ We will use the four preliminary graphs $B_1,$ $C_1,$ $D_1$ and $D_2$ in Figure 3.
\begin{figure}[h]
\centering
\includegraphics[width=0.8\textwidth]{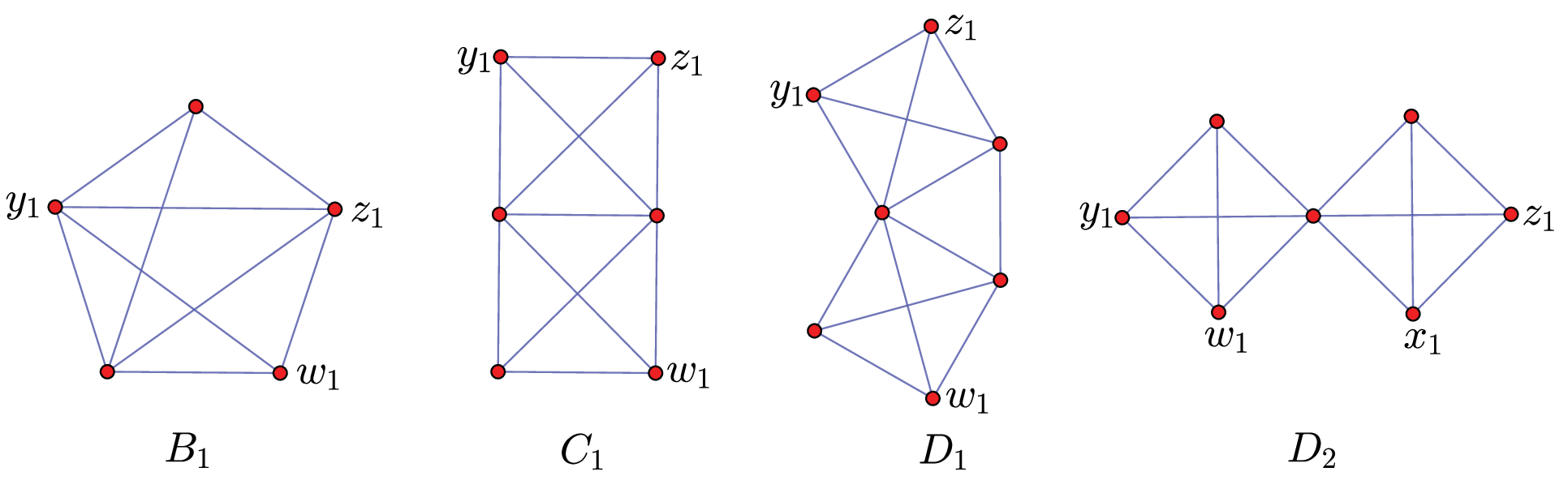}
\caption{$B_1,$ $C_1,$ $D_1$ and $D_2$}
\end{figure}

Given a positive integer $k,$ let $A_i$ be a graph isomorphic to $K_4$ with vertex set $V(A_i)=\{x_i, y_i, z_i, w_i\}$ for $i=1,\ldots, 2k+1.$

Suppose $n=8k.$ Let $G_n$ be the graph obtained from the disjoint union $A_1+A_2+\cdots +A_{2k}$ by adding edges $z_iy_{i+1},$ $1\le i\le 2k$ where $y_{2k+1}=y_1$ and $w_jw_{k+j},$
$1\le j\le k.$

Suppose $n=8k+1.$ Replace $A_1$ by $B_1$ and then construct $G_n$ as above.

Suppose $n=8k+2.$ Replace $A_1$ by $C_1$ and then construct $G_n$ as in the case $n=8k$.

Suppose $n=8k+3.$ Replace $A_1$ by $D_1$ and then construct $G_n$ as in the case $n=8k$.

Suppose $n=8k+4.$ Let $G_n$ be the graph obtained from the disjoint union $A_1+A_2+\cdots +A_{2k+1}$ by adding edges $z_iy_{i+1},$ $1\le i\le 2k+1$ where $y_{2k+2}=y_1$ and adding edges $w_1w_{k+1,}$ $w_jw_{k+j+1},$ $1\le j\le k.$

Suppose $n=8k+5.$ Replace $A_1$ by $B_1$ and then construct $G_n$ as in the case $n=8k+4$.

Suppose $n=8k+6.$ Replace $A_1$ by $C_1$ and then construct $G_n$ as in the case $n=8k+4$.

Suppose $n=8k+7.$ Let $G_n$ be the graph obtained from the disjoint union $D_2+A_2+\cdots +A_{2k+1}$ by adding edges $z_iy_{i+1},$ $1\le i\le 2k+1$ where $y_{2k+2}=y_1$ and adding edges $w_1w_{k+1,}$ $x_1w_{k+2},$ and $w_jw_{k+j+1},$ $2\le j\le k.$

Finally, it is easy to verify that $G_n$ is $3$-connected, locally nonforesty and has size $f(n).$ This completes the proof. \hfill$\Box$

{\bf Remark.} Denote $b(n)=7(n-1)/3,$ the lower bound in Conjecture 1 and let $f(n)$ be defined as in Theorem 1. Clearly, for any $n\ge 8,$
$$
b(n)>2n\ge f(n).
$$
Their difference can be large. For example, $b(100)-f(100)=43.$ Thus Conjecture 1 does not hold.

\vskip 5mm
{\bf Acknowledgement.} This research  was supported by the NSFC grant 12271170 and Science and Technology Commission of Shanghai Municipality
 grant 22DZ2229014.

\end{document}